\documentclass{llncs}
\usepackage{makeidx}  
\usepackage{setspace}
\usepackage{amssymb}
\usepackage{amsmath}
\usepackage{amsfonts}
\usepackage{amssymb}
\usepackage{epsfig}
\usepackage{graphicx}
\usepackage{subfigure}
\usepackage[nosepfour,warning,np,debug,autolanguage]{numprint}
\begin{document}
\pagestyle{headings}
\frontmatter          
\pagestyle{headings}  
\title{Which Meshes are Better Conditioned Adaptive, Uniform, Locally Refined or Locally Adjusted ?}
\titlerunning{Better Conditioned Meshes}  
%
\author{Sanjay Kumar Khattri and Gunnar Fladmark}
%
\authorrunning{S. K. Khattri et al.}   
%
%
\institute{Department of Mathematics, University of Bergen, Norway\\
\email{\{sanjay,\,\,Gunnar.Fladmark\}@mi.uib.no} \\ 
\texttt{{http://www.mi.uib.no/$\sim$sanjay}}
}

\maketitle              
\begin{abstract}
Adaptive, locally refined and locally adjusted meshes are preferred over uniform meshes for capturing singular or localised solutions. Roughly speaking, for a given degree of freedom a solution associated with adaptive, locally refined and locally adjusted meshes is more accurate than the solution given by uniform meshes. In this work, we answer the question which meshes are better conditioned. We found, for approximately same degree of freedom (same size of matrix), it is easier to solve a system of equations associated with an adaptive mesh.
\end{abstract}
\pagestyle{empty} 
\section{Introduction}
\begin{figure*} 
\begin{minipage}[t]{5.0cm}
\includegraphics[scale=0.25]{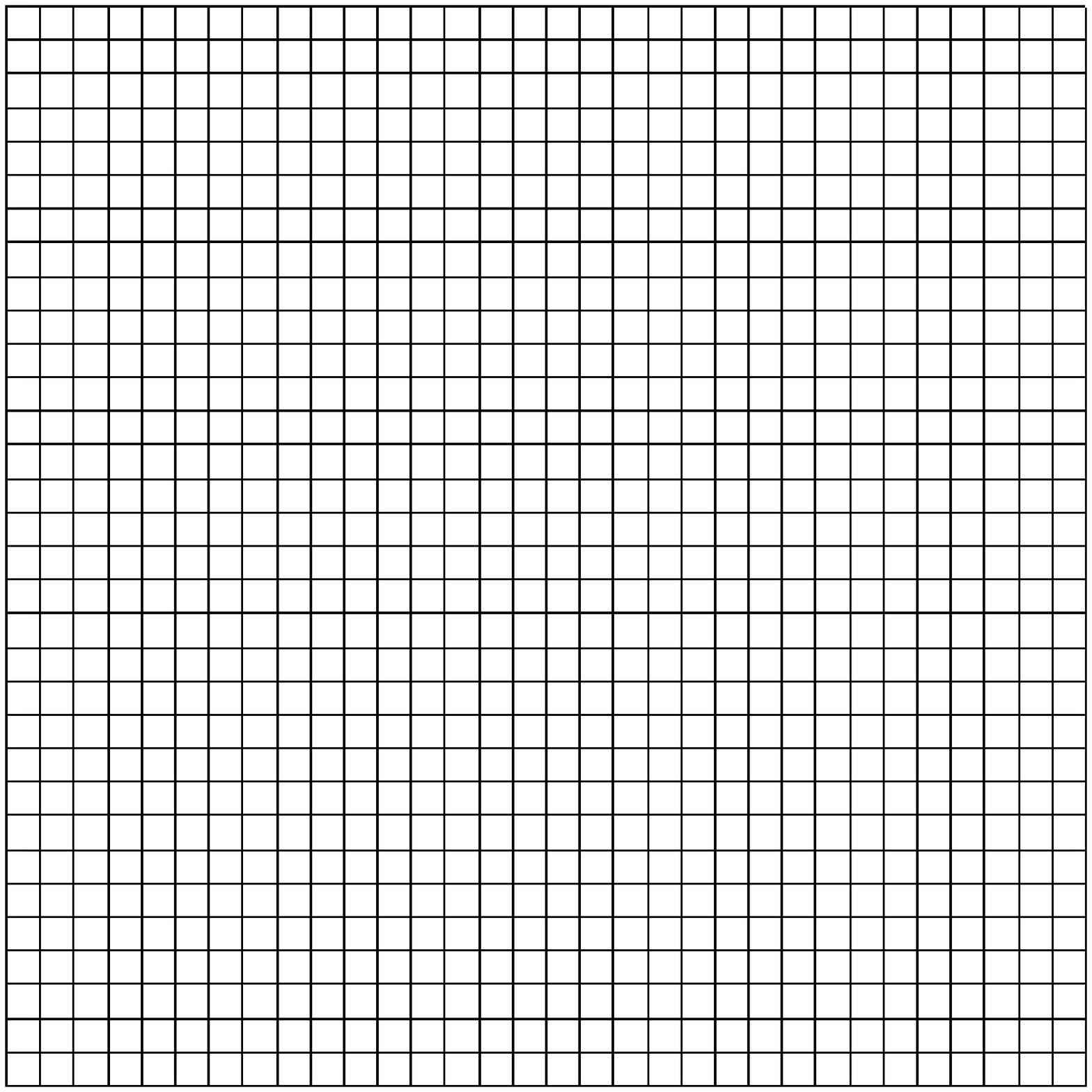} 
\caption{Uniform mesh.} 
\label{fig:uniform_}
\end{minipage}
\hfill
\begin{minipage}[t]{5.0cm}
\includegraphics[scale=0.25]{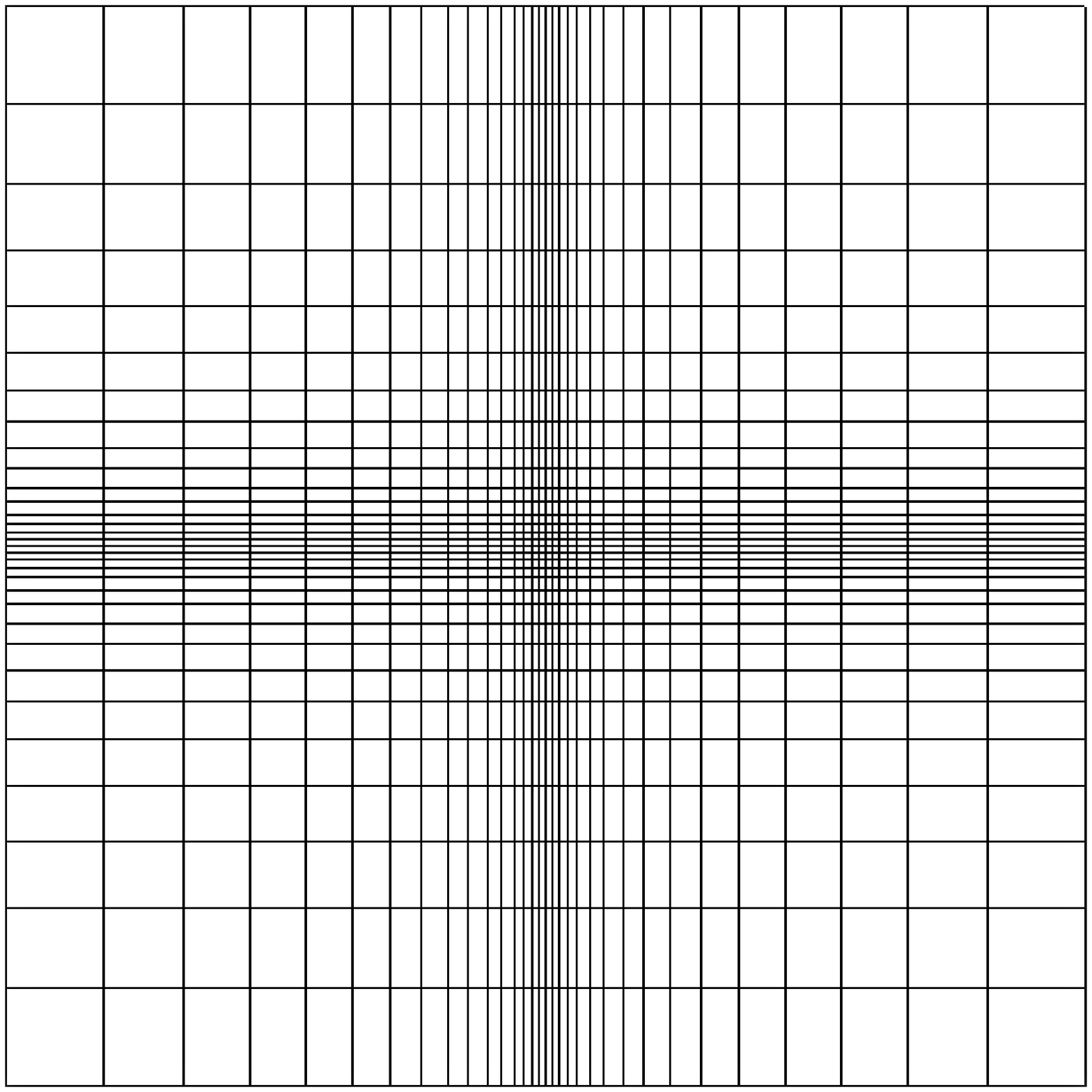}
\caption{Localised mesh.}
\label{fig:localised_}
\end{minipage}
\hfill \\
\begin{minipage}[t]{5.0cm}
\includegraphics[scale=0.25]{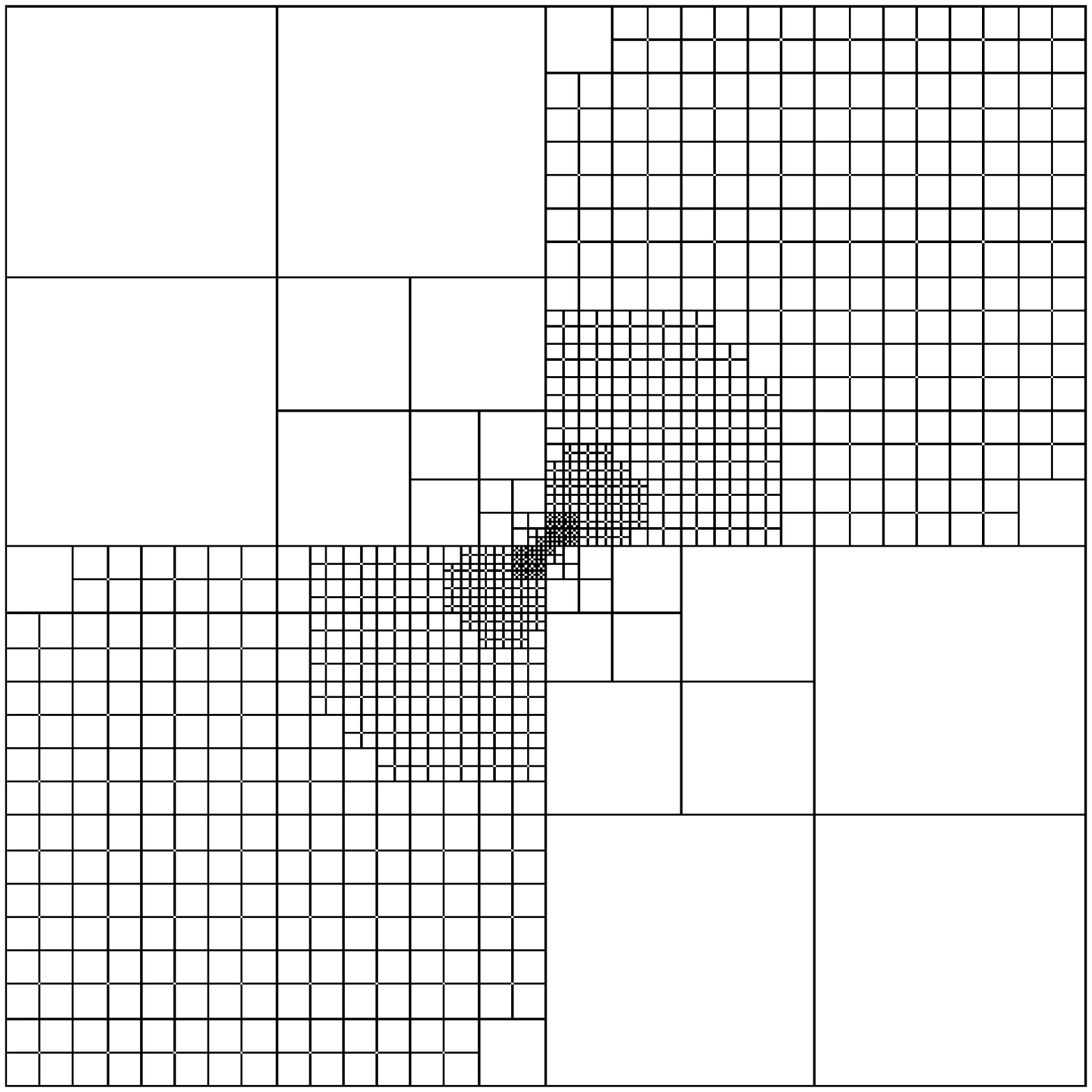} 
\caption{Adaptive mesh.}
\label{fig:adaptive_}
\end{minipage}
\hfill 
\begin{minipage}[t]{5.0cm}
\includegraphics[scale=0.250]{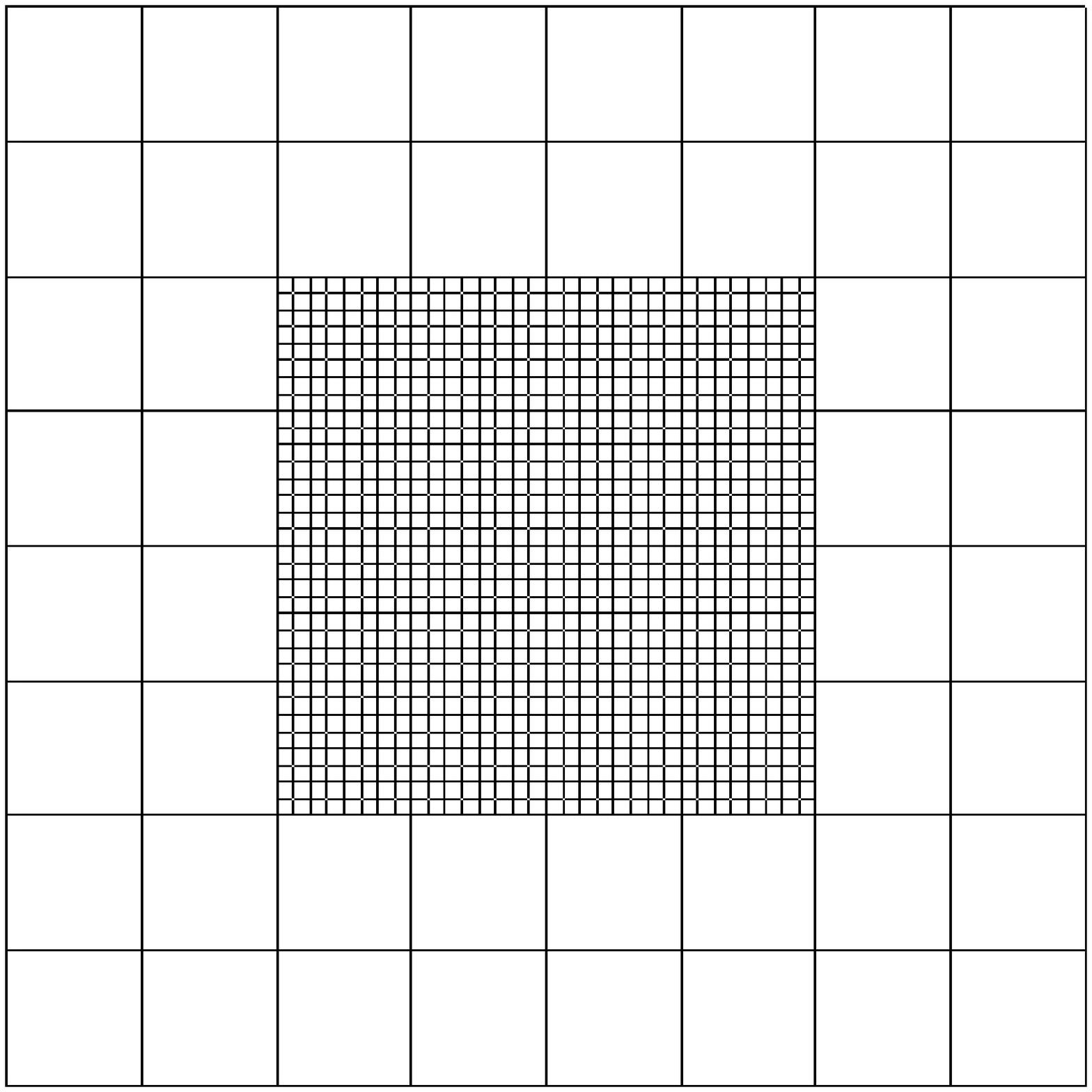}
\caption{Locally refined mesh.}
\label{fig:locallyrefined_}
\end{minipage}
\label{meshes_}
\end{figure*}

Uniform, locally adjusted, adaptive and locally refined meshes are shown in Figures \ref{fig:uniform_}, \ref{fig:localised_}, \ref{fig:adaptive_} and \ref{fig:locallyrefined_}, respectively. Here, for each mesh the number of cells (or degree of freedom) are approximately \numprint{1024} ($2^5\times2^5$). Let us consider the steady state pressure equation of a single phase flow in a porous medium $\Omega$ \cite{Ivar1}
\begin{equation}
-\,\text{div}\,(\mathrm{K}\,\text{grad}\,p) = f \qquad \text{in} \quad \Omega\quad\text{and}\quad
p(x,y) = p^\text{D}             \qquad  \text{on} \quad  \partial{\Omega_\text{D}}\enspace .
\label{elliptic1}
\end{equation}
Here, $\Omega$ is a polyhedral domain in $\mathbb{R}^2$, the source function $f$ is assumed to be in $L^2(\Omega)$ and the diagonal tensor coefficient ${\mathrm{K}{(x,y)}}$ is positive definite and piecewise constant. $\mathrm{K}$ (permeability) is allowed to be discontinuous in space. 
We are discretizating the equation \eqref{elliptic1} on the meshes (see Figures \ref{fig:uniform_}, \ref{fig:localised_}, \ref{fig:adaptive_} and \ref{fig:locallyrefined_}) by the method of Finite Volumes \cite{Ivar1,ewing_3,suli_1,tpfa_00}. For discretization of the problem \eqref{elliptic1} on uniform and localised meshes (see the Figures \ref{fig:uniform_} and \ref{fig:localised_}), we refer to the References  \cite{Ivar1,suli_1,tpfa_00}. Discretization of the equation \eqref{elliptic1} on adaptive and locally refined meshes is given in the following References \cite{ewing_3,sanjay_art}. Finite Volume discretization of the problem \eqref{elliptic1} on a mesh results in a matrix system $\boldsymbol{A}\,\mathbf{p}_{{h}}=\mathbf{b}$. Here, $\boldsymbol{A}$ is symmetric positive definite matrix associated with a mesh.  

\begin{center}
\begin{figure}
\begin{minipage}[b]{0.5\linewidth} 
\begin{center}
\includegraphics[scale=0.50]{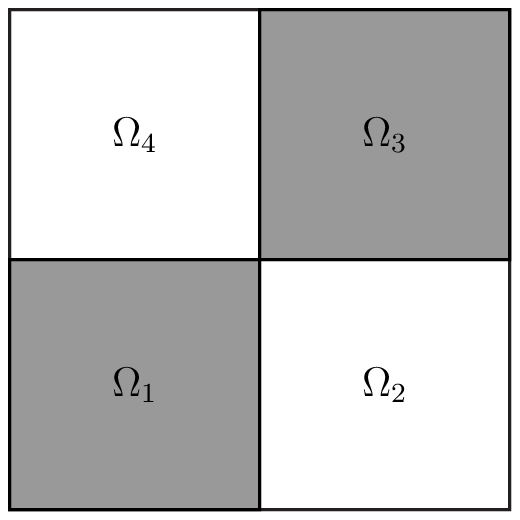}\label{fig:domain_}
\caption{Domain.}
\end{center}
\end{minipage}
\hfill
\begin{minipage}[b]{0.5\linewidth}
\begin{center}
\includegraphics[scale=0.2500]{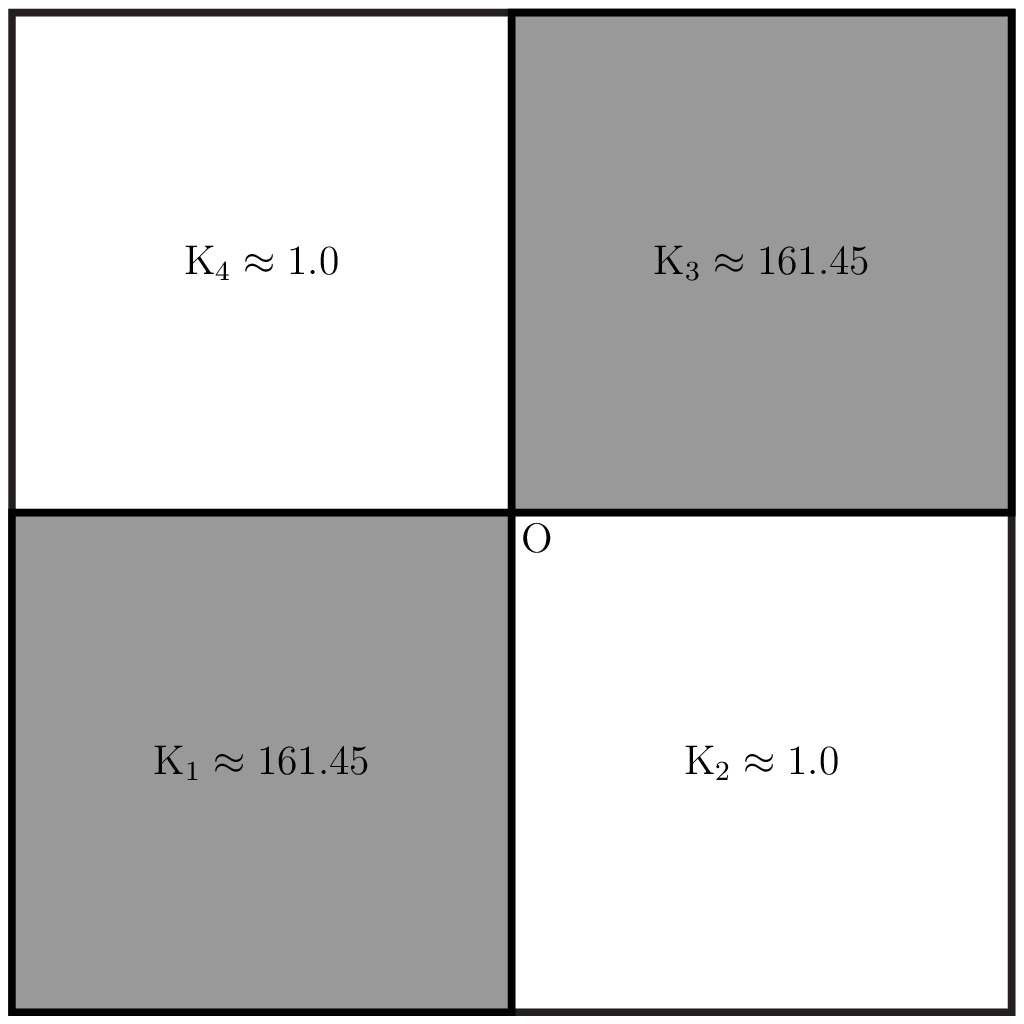}
\label{fig:perm_100_dist}
\caption{Permeability distribution.}
\end{center}
\end{minipage}
\end{figure}
\end{center}

\begin{figure}
\begin{center}
\includegraphics[scale=0.500]{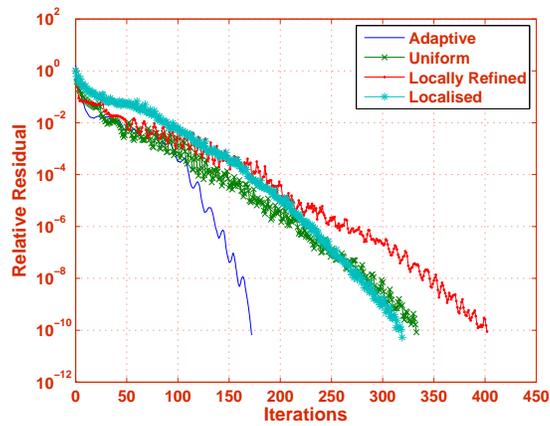} 
\caption{Convergence curves for the matrix system formed on the different meshes.}
\label{fig:comparison}
\end{center}
\end{figure}

Let us define a problem to be solved on the four meshes. Let the domain be $\Omega = [-1,1]\times[-1,1]$ (see Figure \ref{fig:domain_}). It is divided into four sub-domains according to the permeability $\mathrm{K}$ (see the Figures \ref{fig:domain_} and \ref{fig:perm_100_dist}). The permeability $\mathrm{K}$ is a positive constant in each of the sub-domains and is discontinuous across the surfaces of sub-domains. Let the permeability in the sub-domain $\Omega_i$ be $\mathrm{K}_i$. Assuming that $\mathrm{K}_1 = \mathrm{K}_3 = R$ and $\mathrm{K}_2 = \mathrm{K}_4 = 1.0$. $\mathrm{K}_1$, $\mathrm{K}_2$, $\mathrm{K}_3$ and $\mathrm{K}_4$ refers to the permeabilities in the subdomains $\Omega_1$, $\Omega_2$, $\Omega_3$ and $\Omega_4$, respectively. The parameter $R$ is given below. Let the exact solution in the polar form be \cite{sanjay_art} 
\begin{equation}
\label{solu_examp1}
{ p(r,\theta) =r^\gamma \,\eta(\theta)}\enspace,
\end{equation}
where the parameter $\gamma$ denotes the singularity in the solution \cite{sanjay_art} and it depends on the permeability distribution in the domain (see Figure \ref{fig:perm_100_dist} for the permeability for the singularity $\gamma = 0.1$). $\eta(\theta)$ is given as
\begin{equation}
\label{eta_theta_1}
\eta(\theta) = 
\begin{cases}
\cos[(\pi/2-\sigma)\gamma]\, \cos[(\theta-\pi/2+\rho)\gamma]\enspace,  \quad &\theta\in[0,\pi/2]\enspace,\\
\cos(\rho \gamma)\, \cos[(\theta-\pi+\sigma)\gamma]\enspace,           \quad &\theta\in[\pi/2,\pi]\enspace,\\
\cos(\sigma \gamma)\, \cos[(\theta-\pi-\rho)\gamma]\enspace,           \quad &\theta\in[\pi,3\pi/2]\enspace,\\
\cos[(\pi/2-\rho)\gamma]\, \cos[(\theta-3\pi/2-\sigma)\gamma]\enspace, \quad &\theta\in[3\pi/2,2\pi]\enspace.
\end{cases}
\end{equation}
It can be shown that solution $p$ (given by equation \eqref{solu_examp1}) barely belongs in the fractional Sobolev space ${\mathbf{H}^{1+\kappa}(\Omega)}$ with $\kappa < \gamma$ (cf. \cite{strang}). 
For the singularity $\gamma=0.1$, the parameters are
$$R \approx 161.4476\enspace, \quad \rho \approx 0.7854 \quad \text{and}\quad  \sigma \approx  -14.9225\enspace.$$
We solve the problem \eqref{elliptic1} on the four meshes. The exact solution is given by the equation \eqref{solu_examp1}. We enforce the solution inside the domain by the Dirichlet boundary condition and the source term. For solving discrete system of equations formed on the meshes, we use the Conjugate Gradient (CG) solver (see \cite{tech_rep}). Table \ref{table:eig_values} presents eigenvalues and condition numbers of the matrix systems associated with the different meshes. Note that in this table, the largest eigenvalue on all four meshes is approximately same. However, the smallest eigenvalue associated with the adaptive mesh is greater than the smallest eigenvalues associated with other three meshes.  When solving the Symmetric Positive Definite (SPD) linear system $\boldsymbol{A}\,\mathbf{p}_h=\mathbf{b}$ with the CG, the smallest eigenvalues of the matrix slowes down the convergence (cf. \cite{tech_rep}). Several techniques have been proposed in the literature to remove bad effect of the smallest eigenvalue (see \cite[and references therein]{small_eig_00,tech_rep}). Convergence of the CG solver for these the four systems are shown in the Figure \ref{fig:comparison}. It is clear from the Table \ref{table:eig_values} and the Figure \ref{fig:comparison} that it is easier to solve a matrix system associated with an adaptive mesh than to solve systems associated with uniform, localised and locally refined meshes.
\begin{table}[ht]
\setlength{\extrarowheight}{4pt}
\caption{Eigenvalues and condition numbers of different matrix systems.}
\centering
\begin{tabular}{|c!{\vrule width 3pt} c| c| c |}
\hline\hline
Mesh & Smallest eig. & Largest eig. & Cond. Num.\\
\hline
Adaptive & $\mathbf{4.15\times{10}^{-1}}$& $1.28\times{10}^{3}$ &  $\mathbf{3.10\times{10}^{3}}$\\
\hline
Localised &$5.50\times{10}^{-2}$ &$0.78\times{10}^3$ & $1.42\times{10}^4$ \\
\hline
Uniform &$7.62\times{10}^{-2}$ & $1.28\times{10}^{3}$ & $1.69\times{10}^{4}$ \\
\hline
Locally Refined & $3.94\times{10}^{-2}$ & $1.28\times{10}^{3}$&  $3.25\times{10}^{4}$ \\
\hline
\hline
\end{tabular}
\label{table:eig_values}
\end{table}
\section{Conclusions}
\label{sec:conclusion}
We have shown that it is easier to solve a matrix system associated with an adaptive mesh than solving systems associated with uniform, localised and locally refined meshes. The adaptive mesh is generated by equal distribution of the fluxes over all the cells in the mesh. Why do equal distribution of fluxes is create meshes which are better conditioned ? Or, why do equal distribution of fluxes remove bad effect of small eigenvalue ? Answers to these questions can help in designing new preconditioners or improving existing preconditioners.
%
\bibliographystyle{splncs}

\end{document}